\documentclass[reqno,12pt]{article}

\usepackage{a4wide}
\usepackage{amsmath,stmaryrd} 
\usepackage{amssymb}
\usepackage{amsthm}
\usepackage[utf8]{inputenc} 
\usepackage{graphicx} 
\usepackage{xcolor}
\usepackage[normalem]{ulem}

\usepackage{url}

\numberwithin{equation}{section}

\newtheorem{theo}{Theorem}

\theoremstyle{remark}

\newtheorem*{Remark*}{Remark}
\newtheorem*{Remarks*}{Remarks}

\makeatletter
\newcommand*{\house}[1]{%
 \mathord{%
 \mathpalette\@house{#1}%
 }%
}
\newcommand*{\@house}[2]{%
 \dimen@=\fontdimen8 %
 \ifx#1\scriptscriptstyle\scriptscriptfont
 \else\ifx#1\scriptstyle\scriptfont
 \else\textfont\fi\fi
 3 %
 \sbox0{%
 $#1%
 \vrule width\dimen@\relax
 \overline{%
 \kern2\dimen@
 \begingroup 
 #2%
 \endgroup
 \kern2\dimen@
 }%
 \vrule width\dimen@\relax
 \mathsurround=1.5\dimen@ 
 $%
 }%
 \ht0=\dimexpr\ht0-\dimen@\relax
 \dp0=\dimexpr\dp0+2\dimen@\relax
 \vbox{%
 \kern\dimen@ 
 \copy0 %
 }%
}

\newcommand{\Q}{\mathbb{Q}}

\newcommand{\K}{\mathbb{K}}

\newcommand{\Qbar}{\overline{\mathbb Q}}

\newcommand{\popt}{\lambda}

\title{Transcendence measure for the values of the logarithmic derivative of the Bessel function $J_0$ at algebraic arguments}

\author{S. Fischler and T. Rivoal}

\date\today

\begin{document}

\maketitle

\begin{abstract} Let $P\in \mathbb Z[X]\setminus\{0\}$ be of degree $\delta\ge 1$ and usual height $H\ge 1$, and let $\alpha\in \Qbar^*$ be of degree $d\ge 2$. As a consequence of general result due to Lang and Galochkin, we have the following transcendence measure:  for any $\varepsilon>0$, there exists $c>0$ such that $\vert P(J_0'(\alpha)/J_0(\alpha))\vert>c/H^{4d^2\delta+\varepsilon}$ where $J_0$ is the Bessel function. In this paper, we prove that the exponent $4d^2\delta$ can be replaced by a smaller (explicit) quantity $\mu(d,\delta)\le 4d^2\delta-2d\delta-1$. A similar improvement holds more generally for the logarithmic derivative of any $E$-function $f$ of differential order 2 and such that $f$ and $f'$ are homogeneously algebraically independent.   Our method rests upon the optimization of the size of a determinant that appears naturally in the classical Siegel-Shidlovskii method, following the steps of our previous improvement of the transcendence measure of the value $e^\alpha$ for any $\alpha\in \Qbar^*$.  
 \end{abstract}
 
\section{Introduction}

The main result of \cite{firiPMB} is a new transcendence measure for the values of the exponential function at algebraic arguments; the improvement over previous measures came from the optimization of the size of a certain determinant, that had already appeared in the work of Mahler \cite{mahler} and Zheng \cite{zheng}. The goal of this paper is to develop further this optimization method to compute new transcendence measures for values at algebraic arguments of $E$-functions of differential order $\le 2$, or of their logarithmic derivatives.  We embed the field of algebraic numbers $\Qbar$ into $\mathbb C$. We recall the original definition of $E$-functions by Siegel \cite{siegel}: a power series $f(x)=\sum_{n=0}^\infty a_nx^n/n!\in \Qbar[[x]]$ is an $E$-function if $f$ satifies a linear differential equation over $\Qbar(x)$, and for all $\varepsilon>0$, there exists $c(\varepsilon)>0$ such that $\house{a_n}\le c(\varepsilon) n!^{\varepsilon}$ and $\textup{den}(a_0, a_1, \ldots, a_n)\le c(\varepsilon) n!^{\varepsilon}$ for all $n\ge 0$. If we replace these two upper bounds by $C^{n+1}$ for some $C>1$, then this defines the {\em a priori} smaller class of $E$-functions in the strict sense, though it is conjectured that both classes coincide (\cite[p.~715]{andre}). Our results hold for $E$-functions in Siegel's sense.

The hypergeometric nature of $E$-functions of order $\le 2$ is well understood. Consider the confluent hypergeometric function ${}_1F_1[\alpha;\beta;x]:=\sum_{n=0}^\infty (\alpha)_nx^n/(n!(\beta)_n)$, $\alpha \in \mathbb C, \beta\in \mathbb C\setminus \mathbb Z_{\le 0}$, where $(a)_n:=a(a+1)\cdots (a+n-1)$. Shidlovskii mentioned in \cite[p. 184]{shid} that $E$-functions solutions of a homogeneous differential equation over $\Qbar(x)$ of order 1 are of the form $p(x)e^{\omega x}$, with $p\in \Qbar[x]$ and $\omega \in \Qbar$. Gorelov proved in \cite{gor1} that $E$-functions solutions of an inhomogeneous differential equation over $\Qbar(x)$ of order 1 are exactly of the form $a(x){}_1F_1[1;\beta;\lambda x]+b(x)$, with $a,b\in\Qbar[x,x^{-1}]$, $\beta \in \{1\} \cup (\mathbb{Q}\setminus\mathbb{Z})$, $\lambda\in \Qbar$. In \cite{gor2}, he then considered $E$-functions solutions of homogeneous differential equations over $\Qbar(x)$ of order 2: if the equation is reducible over $\Qbar(x)$, then these $E$-functions are of the form  $a(x)e^{\mu x}{}_1F_1[\alpha;\beta; \lambda x]+b(x)e^{\mu x}$ with $a,b\in\Qbar[x,x^{-1}]$, $\lambda, \mu\in \Qbar$, $\beta \in \{1\} \cup (\mathbb Q\setminus \mathbb Z_{\le 0})$ such that $\alpha-\beta\in \mathbb Z$, while if the equation is irreducible over $\Qbar(x)$ they are of the form $a(x)e^{\mu x}{}_1F_1[\alpha;\beta; \lambda x]+b(x)e^{\mu x}{}_1F_1'[\alpha;\beta; \lambda x]$ with $a,b\in\Qbar(x)$, $\lambda\in \Qbar^*$, $\mu\in \Qbar$, $\alpha, \beta \in \mathbb Q\setminus \mathbb Z_{\le 0}$ such that $\alpha-\beta\in \mathbb Z$. See \cite{rivroq1, rivroq2} for different proofs of these results.

Given a transcendental number $\xi\in \mathbb C$, a {\em transcendence measure} of $\xi$ is a non-trivial lower bound of $\vert P(\xi)\vert$ in terms of the height and degree of $P\in\Qbar[X]$.  For $\alpha\in \mathbb K$ where $\mathbb K$ is a number field, the house of $\alpha$ is $\house{\alpha}:=\max_{\sigma}\vert \sigma(\alpha)\vert$ where $\sigma$ runs through all embeddings of $\mathbb K$ into $\mathbb C$; in other words, $\house{\alpha}$ is the maximum of the moduli of $\alpha$ and of all its Galois conjugates over $\mathbb Q$. Given $P(X)=\sum_{j=0}^d a_jX^j \in \Qbar[X]$, we set $H(P):=\max_{j=0, \ldots, d}\house{a_j}$ its usual height. We define the usual height $H(\alpha)$ of  $\alpha\in \Qbar$ as $H(Q)$ where $Q\in \mathbb Z[X]\setminus\{0\}$ is the minimal polynomial of $\alpha$ over $\mathbb Q$ (normalized so that its coefficients are coprime and the leading coefficient is positive). We also let $\deg(\alpha):=\deg(Q)$ and $\textup{den}(\alpha)\ge 1$ be the denominator~of~$\alpha$.

Let $N, d, \delta$ be integers such that $d, \delta,N\ge 1$ and $N\ge \delta d$ and let us define
\begin{equation} \label{eq:defpsi}
\psi(d,\delta, N):=\frac{\delta d^2(N-\delta+1)}{N-\delta d+1}+d(N-\delta+1)-1.
\end{equation}
The parameter $N$ in the definition of $\psi(d,\delta,N)$ is accessory but appears naturally in the proofs of Theorems \ref{theo:1} and \ref{theo:2} stated  below.~(\footnote{The parameters we denote here by $N, N_1$ and $N_2$ are denoted in \cite{firiPMB} by $p, p_1$ and $p_2$ respectively.})  If $d=1$, note that $\psi(1,\delta, N)=N$ is independent of $\delta$, but we still require that $N\ge \delta$.  For $d\ge 2$, we define $N_1:=\delta d-1+\lfloor \delta\sqrt{d^2-d}\rfloor$ and $N_2:=N_1+1$, which are both $\ge \delta d$, and we then set 
$
\popt :=
 N_1$ if $\psi(d,\delta,N_1)\le \psi(d,\delta, N_2)$ and $\popt :=
N_2$ if  $\psi(d,\delta,N_2)<\psi(d,\delta, N_1).
$ 
When $d=1$, we set $\popt:=\delta$. We shall prove that $\popt$, which depends on $d$ and $\delta$, minimizes the function $N\mapsto \psi(d,\delta,N)$ amongst all integer values of  $N$ such that $N\ge \delta d$ (this is obvious if $d=1$ by definition). For simplicity, we define $\mu(d,\delta):=\psi(d,\delta, \lambda)$. For $d=1$ and all $\delta\ge 1$, $\mu(1,\delta)=\delta$, while for all $d\ge2$ and all $\delta \ge1$, we have 
\begin{multline*}
\big(2d^2+2d\sqrt{d^2-d}-d\big)\delta-1
\le  
\mu(d,\delta)
\le \big(2d^2+2d\sqrt{d^2-d}-d\big)\delta-1+\frac{d}{\delta \sqrt{d^2-d}-1}.
\end{multline*}  

We shall prove the following two theorems. They are not new when $d=1$ (and thus $\mu(d,\delta)=\delta$) as they can be found in \cite[Chapter 11]{shid}, and this case is recorded here for completeness. On the other hand, when $d\ge 2$, the best known exponent of $H$ in  \eqref{eq:transcmeasure} and \eqref{eq:transcmeasure2} below was $4 d^2\delta$ (the  bound obtained by Lang-Galochkin, see~\cite[p. 238, Theorem 5.29 and remarks]{feldnest} or \cite[p.~383, Theorem~1]{shid}): our results are better because $\mu(d,\delta)\le 4d^2\delta-2d\delta-1$ for all $\delta\ge 1$, and even $\mu(d,\delta)\le 4d^2\delta-(2d+\frac14)\delta$ when $\delta\ge 4$. We refer the reader to the introduction of \cite{firiPMB} for more details, in particular explicit expressions  for $\delta\in\{1, 2,3,4\}$ of $\mu(d,\delta)\in \mathbb Q(d)$ valid for all $d\ge 2$. 
It is also proved in \cite{firiPMB} that for fixed $d$, $\mu(d,\delta)$ is an increasing function of $\delta$, so that the smallest possible exponent of $H$ in \eqref{eq:transcmeasure} and \eqref{eq:transcmeasure2} below is when $\deg(P)=\delta$.
\begin{theo} \label{theo:1} Let $\K$ be a number field and $f\in \K[[x]]$ be an $E$-function solution of a differential equation $Ay''=By'+Cy$ with $A,B,C\in \K[x]$, $A\neq 0$. Let us assume that $f$ and $f'$ are homogeneously algebraically independent over $\Qbar(x)$. Let $\alpha\in \Qbar$ be such that $\alpha A(\alpha)\neq 0$ and set $d:=[\mathbb K(\alpha):\mathbb Q]\ge 1$.

Then for any $\varepsilon>0$ and  any integer $\delta \ge 1$, there exists a constant $c=c(\varepsilon, \alpha, \delta, \mathbb K)>0$ such that for all $H\ge 1$, we have 
\begin{equation}\label{eq:transcmeasure}
  \left\vert P\big(f'(\alpha)/f(\alpha)\big)\right\vert  > \frac{c}{H^{
 \mu(d,\delta)
  +\varepsilon}}
\end{equation}
for every polynomial $P\in \mathbb Z[X]\setminus\{0\}$ of degree $\le \delta$ and height $H(P)\le H$.
\end{theo}
Theorem \ref{theo:1} applies in particular to the case where $f$ is the Bessel function $J_0(x):=\sum_{n=0}^\infty (-1)^n (x/2)^{2n}/n!^2=e^{-ix}{}_1F_1[1/2;1;2ix]$ (solution of $xy''+ y'+xy=0$)  for any $\alpha\in \Qbar^*$, with $\K=\mathbb Q$ and $d=\deg(\alpha)$. We also obtain the optimal transcendence exponent $\mu(1,\delta)=\delta$ for the continued fraction $[1;2,3,4,5,\ldots]=\frac{iJ_0'(2i)}{J_0(2i)}=(\frac{(J_0(2ix))'}{2J_0(2ix)})_{\vert x=1}$ (see \cite[p. 218]{shid}) with $\K=\mathbb Q$, $\alpha=1$ and $d=1$, because $J_0(2ix)$ is an $E$-function in $\mathbb Q[[x]]$. This result is not new however, as it is a consequence of Shidlovskii's $\mathbb Q$-linear independence measure of values at rational arguments of $\mathbb Q(x)$-linearly independent $E$-functions in $\mathbb Q[[x]]$ that form a solution of a differential system (see \cite[p. 357, Theorem 1]{shid}): we apply this result to the $\mathbb Q(x)$-independent functions $f^jf'^{\delta-j}$ with $f(x)=J_0(2i x)$.

Theorem \ref{theo:1} can be slightly generalized as follows: the transcendence measure \eqref{eq:transcmeasure} holds for $f_2(\alpha)/f_1(\alpha)$ where $f_1, f_2$ are homogeneously algebraically independent $E$-functions in $\K[[x]]$ such that $Y:={}^t(f_1, f_2)$ is solution of a differential system $Y'=UY$ with $A\in \K[x]$ a common denominator of the entries of $U\in M_2(\K(x))$, and $\alpha\in \Qbar$ such that $\alpha A(\alpha)\neq 0$. The changes to the proof given in \S\ref{sec:determinant} are minor: one just has to apply the Siegel-Shidlovskii method recalled in \S\ref{sec:ssmethod} to the functions $f_1^{N-j}f_2^j$, $j=0,\ldots, N$, instead of $f^{N-j}f'^j$, $j=0,\ldots, N$.  

\medskip

The assumption that $f$ and $f'$ are homogeneously algebraically independent over $\Qbar(x)$ in Theorem~\ref{theo:1} is of course important to ensure the transcendence of $f'(\alpha)/f(\alpha)$ by the Siegel-Shidlovskii theorem. If this assumption is not satisfied, we now want to know which homogeneous algebraic relation $f$ and $f'$ can satisfy when a transcendental $E$-function $f$ is of differential order $\le 2$. (If $f$ is not transcendental, it is a polynomial, hence of minimal order 1, and we are not interested in this case.) Since {\em a fortiori} $f$ and $f'$ are algebraically dependent over $\Qbar(x)$, Theorem 3 of \cite{rivroq1} applies~(\footnote{Strictly speaking, the results in \cite{rivroq1} are proven for $E$-functions in the strict sense. However they hold more generally for $E$-functions in Siegel's sense because \cite[Proposition 1]{rivroq1} holds {\em mutatis mutandis} for these functions by the results proved in \cite{lepetit}.}) and it follows that $f$ is either of the form $a(x)e^{\alpha x}+b(x)e^{\beta x}$ ($a,b\in \Qbar[x,x^{-1}]$ and $\alpha,\beta\in \Qbar$) or of the form $a(x){}_1F_1[1;\gamma; \alpha x]+b(x)$ ($a,b \in \Qbar[x,x^{-1}]$, $\alpha\in \Qbar$, $\gamma\in \mathbb Q \setminus \mathbb Z$). Note that $e^{x}$ and ${}_1F_1[1;\gamma; x]$ are solutions of the (in)homogeneous linear equation  $xy'=(x+1-\gamma)y+\gamma-1$ of order 1 over $\Qbar(x)$ (take $\gamma=1$ for $e^x$). Therefore, the existence of a homogeneous algebraic relation over $\Qbar(x)$ between $f$ and $f'$ implies certain restrictions on $a,b,\alpha, \beta$, and in turn this implies that $f$ and $f'$ are linearly dependent over $\Qbar(x)$. In other words, $f$ satisfies a  homogeneous linear equation of order 1 over $\Qbar(x)$. Our 
second theorem deals with this case, and more generally with inhomogeneous linear equations of order 1 over $\Qbar(x)$.  

\begin{theo} \label{theo:2} Let $\K$ be a number field and $f\in \K[[x]]$ be an $E$-function solution of a differential equation $Ay'=By+C$ with $A,B,C\in \K[x]$, $A\neq 0$. Let us assume that $f$ is transcendental over $\Qbar(x)$. Let $\alpha\in \Qbar$ be such that $\alpha A(\alpha)\neq 0$ and set $d:=[\mathbb K(\alpha):\mathbb Q]\ge 1$.

Then for any $\varepsilon>0$ and  any integer $\delta \ge 1$, there exists a constant $c=c(\varepsilon, \alpha, \delta, \mathbb K)>0$ such that for all $H\ge 1$, we have 
\begin{equation}\label{eq:transcmeasure2}
  \left\vert P\big(f(\alpha)\big)\right\vert  > \frac{c}{H^{
 \mu(d,\delta)
  +\varepsilon}}
\end{equation}
for every polynomial $P\in \mathbb Z[X]\setminus\{0\}$ of degree $\le \delta$ and height $H(P)\le H$.
\end{theo}
If we take $C=0$ in Theorem \ref{theo:2}, then $f$ is necessarily of the form $f(x)=p(x)e^{\omega x}$ with $\omega \in \Qbar$ and $p(x)\in \Qbar[x]$, and the conclusion is contained in that of \cite[Theorem~1]{firiPMB}. Theorem~\ref{theo:2} applies for all $\gamma\in \mathbb Q\setminus \mathbb Z_{\le 0}$ and all $\alpha\in \Qbar^*$ to ${}_1F_1[1;\gamma; x]=\sum_{n=0}^\infty x^n/(\gamma)_n$, and this is essentially the only example.

\medskip

With minor modifications, one easily proves generalizations of Theorems \ref{theo:1} and \ref{theo:2} in which $P\in \mathcal{O}_{\widetilde{\K}}[X]$ for some number field $\widetilde{\K}$. Both results hold {\em mutatis mutandis} with $d=[\mathbb L: \Q]$ where $\mathbb L$ is the compositum of $\K(\alpha)$ and $\widetilde{\K}$. This is the more general version proved in \cite[Theorem 1]{firiPMB} in the case $f=\exp$.

It would be of course very interesting to improve on the general transcendence measure of  Lang-Galochkin for the values of $E$-functions of differential order $\ge 3$ evaluated at their non-singular points $\alpha$. However, it is not clear to us if the determinant method used in this paper could be adapted to this more general situation.  

\bigskip

The structure of this paper is as follows. We first recall the output of the Siegel-Shidlovskii method in \S \ref{sec:ssmethod}. Then we move in \S \ref{sec:determinant} to the proof of Theorem \ref{theo:1}. At last, we explain in \S \ref{sec:determinant2} the modifications needed to prove Theorem \ref{theo:2}.

\section{A quick reminder of the Siegel-Shidlovskii method}\label{sec:ssmethod}

The following setting is common to the proofs of Theorems \ref{theo:1} and \ref{theo:2}.
Let $f_0, \ldots, f_{N}$ be $E$-functions in $\K[[x]]$ such that $Y:={}^t(f_0, \ldots, f_{N})$ is solution of a linear differential system $Y'=UY$ with $U\in M_{N+1}(\K(x))$. 
Let $T$ denote the least common denominator in $\K[x]\setminus\{0\}$ of the entries of $U$. 
Then by \cite[p. 114, Lemma 16]{shid}, we get the following fundamental result (with slightly modified notations) which is at the core the Siegel-Shidlovskii method. Assume that $f_0, \ldots, f_{N}$ are $\Qbar(x)$-linearly independent and let $\alpha\in \Qbar$ be such that $\alpha T(\alpha)\neq 0$. Then for any $\varepsilon \in (0,1)$ and any $n\ge n_0$, there exists a set of $N+1$ linearly independent~(\footnote{in the sense that the matrix $(A_{n,j,k})_{0\le j,k\le N}$, which depends on $\varepsilon$ and $\alpha$, is invertible for all $n\ge n_0$}) linear forms 
$$
R_{n,k}:=\sum_{j=0}^{N} A_{n,k,j} f_j(\alpha), \quad k=0,\ldots, N,
$$
with $A_{n,k,j}\in \mathcal{O}_{\K(\alpha)}$ such that $\house{A_{n,k,j}}=\mathcal{O}(n!^{1+\varepsilon})$ and $\vert R_{n,k}\vert=\mathcal{O}(n!^{-(N-\varepsilon)})$ as $n\to +\infty$. The integer $n_0$ is the notorious {\em Shidlovskii's constant} (independent of $\alpha)$, the value of which is not important for our applications. In the sequel, we shall write $R_k$ and $A_{k,j}$ for simplicity.

\medskip

To prove Theorem \ref{theo:1}, we shall apply this contruction to the $E$-functions $f_j:=f^{N-j}f'^{j}$, for $j=0,\ldots, N$, which are linearly independent over $\Qbar(x)$ because $f$ and $f'$ are assumed to be homogeneously algebraically independent. Moreover ${}^t(f_0,\ldots, f_N)$ is solution of a linear differential system with coefficients in $\K(x)$ whose common denominator is $A$. Indeed using the differential equation $Af''=Bf'+Cf$, we have 
$$
f_0'=Nf'f^{N-1}=Nf_{1}, 
\quad f_N'=Nf''f'^{N-1}=Nf'^{N-1}\left(\frac{B}{A}f'+\frac{C}{A}f\right) = 
N\frac{B}{A}f_N+N\frac{C}{A}f_{N-1}
$$
and for $1\le j\le N-1$ we have 
\begin{align*}
f_j'&=(N-j)f^{N-j-1}f'^{j+1}+jf^{N-j}f'^{j-1}f''
\\
&=(N-j)f_{j+1}+jf^{N-j}f'^{j-1}\left(\frac{B}{A}f'+\frac{C}{A}f\right)
\\
&=(N-j)f_{j+1}+j\frac{B}{A}f_j+j\frac{C}{A}f_{j-1}.
\end{align*}
In other words, 
$$f'_j=(N-j)f_{j+1}+j\frac{B}{A}f_j+j\frac{C}{A}f_{j-1}
$$ for all $j=0, \ldots, N$ with the conventions that $f_{-1}=f_{N+1}=0$.

\medskip

On the contrary, to prove Theorem \ref{theo:2}, we shall apply this contruction to the $E$-functions $f_j:=f^{j}$, for $j=0,\ldots, N$, which are linearly independent over $\Qbar(x)$ because $f$ is assumed to be a transcendental function. Moreover ${}^t(f_0,\ldots, f_N)$ is solution of a linear differential system with coefficients in $\K(x)$ whose common denominator is $A$. Indeed using the differential equation $Af'=Bf+C$, we have 
$$
f_j'=jf^{j-1}f'=jf^{j-1}\left(\frac{B}{A}f+\frac{C}{A}\right)=j\frac{B}{A}f_j+j\frac{C}{A}f_{j-1}
$$ for all $j=0, \ldots, N$ 
with the convention that $f_{-1}=0$.

\section{Proof of Theorem \ref{theo:1}} \label{sec:determinant}

Let $\delta \ge 1$ and $N\ge \delta$. (Later on, we shall even impose that $N\ge d\delta$ to obtain the result.) We set $P(X):=\sum_{k=0}^\delta a_k X^k$ with $(a_0,\ldots, a_\delta)\in {\mathbb Z}^{\delta+1}\setminus \{0\}$ such that $H(P):=
\max\vert a_k \vert\le H$. Let $\alpha\in \Qbar$ such that $\alpha A(\alpha)\neq 0$.  Since $f$ and $f'$ are homogeneously algebraically independent over $\Qbar(x)$, $f(\alpha)$ and $f'(\alpha)$ are homogeneously algebraically independent over $\Qbar$ by the homogeneous version of the  Siegel-Shidlovskii theorem \cite[p. 83]{shid} applied to the vector ${}^t(f,f')$. In particular  $f(\alpha)$ is non-zero, and $f'(\alpha)/f(\alpha)$ is a transcendental number. 

The $N-\delta+1$ vectors  ${}^t(a_0, \ldots, a_\delta, 0,\ldots, 0)$, ${}^t(0, a_0, \ldots, a_\delta, 0,\ldots, 0)$,..., ${}^t(0,\ldots, 0,a_0, \ldots, \break a_{\delta})$ of $\mathbb C^{N+1}$  are $\mathbb C$-linearly independent. Since the matrix $(A_{k,j})_{0 \le j,k \le N}$ (constructed in \S\ref{sec:ssmethod} with $f_j:=f^{N-j}f'^{j}$)  is invertible, we can complete these vectors with $\delta$  
vectors ${}^t(A_{\ell_j,0}, A_{\ell_j,1}, \ldots, A_{\ell_j,N})$ ($j=1, \ldots, \delta$, $\ell_j\in \{0,\ldots, N\}$) to form a basis of $\mathbb C^{N+1}$. Up to renumbering and for simplicity, we assume from now on without loss of generality that $\ell_j=j$ for $j=0,\ldots, \delta-1$.

It follows that the algebraic integer of $\mathbb K(\alpha)$ 
\begin{equation} \label{eqDun}
    D:=\left\vert 
\begin{matrix}
a_0&a_1&\cdots &a_{\delta}&0&\cdots&\cdots &\cdots&\cdots &0
\\
0 &a_0&\cdots &a_{\delta-1} &a_{\delta}&0&\cdots &\cdots&\cdots &0
\\
\vdots &\vdots&\vdots &\vdots&\vdots&\vdots&\vdots&\vdots&\vdots&\vdots
\\
0 &0&\cdots &0&\cdots&\cdots &0&a_0&\cdots & a_{\delta}
\\
A_{0,0}&A_{0,1}&\cdots &\cdots&\cdots&\cdots &\cdots &\cdots &\cdots& A_{0,N}
\\
A_{1,0}&A_{1,1}&\cdots &\cdots &\cdots \cdots &\cdots &\cdots&\cdots&\cdots& A_{1,N}
\\
\vdots &\vdots&\vdots &\vdots&\vdots&\vdots&\vdots&\vdots&\vdots&\vdots
\\
A_{\delta-1,0}&A_{\delta-1,1}&\cdots &\cdots&\cdots \cdots &\cdots&\cdots &\cdots &\cdots& A_{\delta-1,N}
\end{matrix}
\right\vert
\end{equation}
is non-zero. For every embedding $\sigma$ of $\mathbb K(\alpha)$ into $\mathbb C$, we also have $\sigma(D)\neq 0$ so that
$$
\prod_{\sigma} \sigma(D) \in \mathbb Z \setminus \{0\}.
$$
In this product, which is over all such embeddings, we shall distinguish $D$ from the other $\sigma(D)$ with $\sigma\neq id$.~(\footnote{To prove the generalization of both Theorems \ref{theo:1} and \ref{theo:2} mentioned at the end of the introduction, the changes to be made in the proofs are as follows: we let $\mathbb L$ be the compositum of $\K(\alpha)$ and $\widetilde{\K}$, $d$ is equal to $[\mathbb L: \mathbb Q]$, the $a_j$'s are assumed to be in $\widetilde{\K}$, the embeddings $\sigma$ are those of $\mathbb L$ into $\mathbb C$, and the entries $a_j$ in the determinant $\sigma(D)$ must be replaced by $\sigma(a_j)$. })

We define for simplicity $\beta:=f'(\alpha)/f(\alpha)$. Let $L_j:=\sum_{k=0}^\delta a_k \beta^{k+j}$. On the one hand, by linear combinations of columns, we find  
$$
D=\left\vert 
\begin{matrix}
L_0&a_1&\cdots &a_{\delta}&0&\cdots&\cdots&\cdots &\cdots &0
\\
L_1 &a_0&\cdots &a_{\delta-1}&a_{\delta}&0&\cdots&\cdots&\cdots &0
\\
\vdots &\vdots&\vdots &\vdots&\vdots&\vdots&\vdots&\vdots&\vdots&\vdots
\\
L_{N-\delta} &0&\cdots &0&\cdots&\cdots&0&a_0&\cdots & a_{\delta}
\\
f(\alpha)^{-N}R_{0}&A_{0,1}&\cdots &\cdots &\cdots&\cdots&\cdots &\cdots&\cdots& A_{0,N}
\\
f(\alpha)^{-N}R_{1}&A_{1,1}&\cdots &\cdots &\cdots&\cdots&\cdots &\cdots&\cdots& A_{1,N}
\\
\vdots &\vdots&\vdots &\vdots&\vdots&\vdots&\vdots&\vdots&\vdots&\vdots
\\
f(\alpha)^{-N}R_{\delta-1}&A_{\delta-1,1}&\cdots &\cdots &\cdots& \cdots&\cdots &\cdots &\cdots& A_{\delta-1,N}
\end{matrix}
\right\vert.
$$
Expanding this determinant along its first column and since $L_j=\beta^j L_0$, the bounds given in \S\ref{sec:ssmethod} above 
imply that, for all $\varepsilon\in (0,1)$, all $n\ge 0$ and $H \ge 1$, 
\begin{equation} \label{eq:majorationD}
\vert D\vert \le C_1 H^{N-\delta}n!^{\delta(1+\varepsilon)} \vert L_0\vert + C_1 \frac{H^{N-\delta+1}}{n!^{N-\delta-1-\delta\varepsilon}}
\end{equation}
for a constant $C_1>0$ independent of $n$ and $H$.

On the other hand, since 
\begin{equation}
    \label{eqsigmaDun}
\sigma(D)=\left\vert 
\begin{matrix}
a_0&a_1&\cdots &a_{\delta}&0&\cdots&\cdots &\cdots&\cdots &0
\\
0 &a_0&\cdots &a_{\delta-1} &a_{\delta}&0&\cdots &\cdots &\cdots &0
\\
\vdots &\vdots&\vdots &\vdots&\vdots&\vdots&\vdots&\vdots&\vdots&\vdots
\\
0 &0&\cdots &0&\cdots&\cdots&0&a_0&\cdots &a_{\delta}
\\
\sigma(A_{0,0})&\sigma(A_{0,1})&\cdots &\cdots&\cdots &\cdots &\cdots&\cdots &\cdots& \sigma(A_{0,N})
\\
\sigma(A_{1,0})&\sigma(A_{1,1})&\cdots &\cdots &\cdots &\cdots &\cdots&\cdots&\cdots& \sigma(A_{1,N})
\\
\vdots &\vdots&\vdots &\vdots&\vdots&\vdots&\vdots&\vdots&\vdots&\vdots
\\
\sigma(A_{\delta-1,0})&\sigma(A_{\delta-1,1})&\cdots &\cdots&\cdots&\cdots &\cdots &\cdots &\cdots& \sigma(A_{\delta-1,N})
\end{matrix}
\right\vert, 
\end{equation}
we have  for all $\varepsilon\in (0,1)$, all $n\ge 0$ and $H \ge 1$, 
\begin{equation} \label{eq:majorationsigmaD}
\vert \sigma(D)\vert \le C_1H^{N-\delta+1}n!^{\delta(1+\varepsilon)}
\end{equation}
where the constant $C_1>0$ can be taken the same as before (up to  increasing it if necessary). 

Therefore, from $\vert D\vert \prod_{\sigma\neq id} \vert \sigma(D)\vert \ge 1$ and with $d:=[\mathbb K(\alpha):\mathbb Q]$, we deduce using 
\eqref{eq:majorationD} and \eqref{eq:majorationsigmaD}
that 
$$
\frac{1}{\big(C_1 H^{N-\delta+1}n!^{\delta(1+\varepsilon)}\big)^{d-1}} \le \vert D\vert \le C_1 H^{N-\delta}n!^{\delta(1+\varepsilon)} \vert L_0\vert + \frac{C_1H^{N-\delta+1}}{n!^{N-\delta+1-\delta \varepsilon}}.
$$
Hence, 
\begin{equation}\label{eq:minL}
\vert L_0\vert \ge \frac{1}{C_1^{d}H^{(N-\delta+1)(d-1)+N-\delta}n!^{\delta d(1+\varepsilon)}}-\frac{H}{n!^{N+1}}=:M_1
\end{equation}
The right-hand side of \eqref{eq:minL} satisfies
\begin{equation}\label{eq:minM}
M_1\ge \frac{1}{2C_1^{d}H^{(N-\delta+1)(d-1)+N-\delta}n!^{\delta d(1+\varepsilon)}}
\end{equation}
provided we can choose $n\ge n_0$ (minimal to get a lower bound as large as possible) such that 
\begin{equation}\label{eq:minN}
    2H^{(N-\delta+1)d}\le n!^{N+1-d\delta(1+\varepsilon)} C_1^{-d}.
\end{equation}
Since $H\ge 1$ is arbitrary, a necessary condition for the existence of such an $n$ in all circumstances is that $N-\delta d+1>d\delta\varepsilon$, {\em i.e.}, that $N\ge \delta d$ because they are integers (assuming as we may from the beginning that $d\delta\varepsilon<1$). We thus now assume that $N\ge \delta d$ and choose $n\ge n_0$ minimal such that 
\eqref{eq:minN} is satisfied. 
Combining \eqref{eq:minL} and \eqref{eq:minM} with this value of $n$, by standard computations (see \cite[p.~359]{shid} or \cite[\S3.1]{firiPMB}), we finally obtain that for all $\varepsilon>0$, there exists a constant $c=c(\varepsilon, \alpha, \delta, \mathbb K)>0$ independent of $H$ such that
$$
\vert L_0 \vert \ge \frac{c}{H^{\psi(d,\delta,N)+\varepsilon}},
$$
where 
\begin{equation*} 
\psi(d,\delta, N):=\frac{\delta d^2(N-\delta+1)}{N-\delta d+1}+d(N-\delta+1)-1
\end{equation*}
is the function defined in \eqref{eq:defpsi} in the introduction. 
It remains to find the minimal possible value of $\psi(d,\delta, N)$ under the assumption that $N\ge \delta d$. We recall that when $d=1$, $\psi(1,\delta, N)=\delta$ for all $N\ge \delta\ge 1$ so that the minimal value of $\psi(d,\delta, N)$ is achieved at $N=\delta$. We now assume that $d\ge 2$ and $\delta\ge 1$ are fixed. Then the minimum of $x\mapsto \psi(d,\delta,x)$ is attained at 
$$
x_0:=\delta d-1+\delta\sqrt{d^2-d},
$$
and the integers $N_1:=\lfloor x_0\rfloor=\delta d-1+\lfloor \delta\sqrt{d^2-d}\rfloor$ and $N_2:=\lfloor x_0\rfloor+1=\delta d+\lfloor \delta\sqrt{d^2-d}\rfloor$ are both admissible to minimize $\psi(d,\delta,x)$ with $x$ an integer (because both are $\ge \delta d$). Therefore defining $\popt$ as either $N_1$ if $\psi(d,\delta,N_1)\le \psi(d,\delta, N_2)$ or $N_2$ if $\psi(d,\delta, N_2)<\psi(d,\delta, N_1)$, we obtain that $$\psi(d,\delta, N)\ge \psi(d,\delta,\popt)=:\mu(d,\delta)$$ for all $N\ge \delta d$.  This completes the proof of Theorem~\ref{theo:1}.

\section{Proof of Theorem \ref{theo:2}} \label{sec:determinant2}

The proof is very similar to that of Theorem \ref{theo:1}. Let $\delta \ge 1$ and $N\ge \delta$. We set $P(X):=\sum_{k=0}^\delta a_k X^k$ with $(a_0,\ldots, a_\delta)\in {\mathbb Z}^{\delta+1}\setminus \{0\}$ such that $H(P):=
\max\vert a_k \vert\le H$. Let $\alpha\in \Qbar$ such that $\alpha A(\alpha)\neq 0$.  Since $f$ is transcendental over $\Qbar(x)$, $f(\alpha)$ is transcendental over $\Qbar$ for all $\alpha\in \Qbar$ such that $\alpha A(\alpha)\neq 0$ by the Siegel-Shidlovskii theorem \cite[p. 123]{shid} applied to the vector ${}^t(1, f)$.

The $N-\delta+1$ vectors  ${}^t(a_0, \ldots, a_\delta, 0,\ldots, 0)$, ${}^t(0, a_0, \ldots, a_\delta, 0,\ldots, 0)$,..., ${}^t(0,\ldots, 0,a_0, \ldots, \break a_{\delta})$ of $\mathbb C^{N+1}$  are $\mathbb C$-linearly independent. Since the matrix $(A_{k,j})_{0 \le j,k \le N}$ (constructed in \S\ref{sec:ssmethod}, but now with  $f_j:=f^{j}$) is invertible, we can complete these vectors with $\delta$ distinct  $\mathbb C$-linearly independent vectors ${}^t(A_{\ell_j,0}, A_{\ell_j,1}, \ldots, A_{\ell_j,N})$ ($j=1, \ldots, \delta$, $\ell_j\in \{0,\ldots, N\}$) to form a basis of $\mathbb C^{N+1}$. Up to renumbering and for simplicity, we assume from now on without loss of generality that $\ell_j=j$ for $j=0,\ldots, \delta-1$.

We define $D$ as in Eq. \eqref{eqDun}; it is again a non-zero algebraic integer of $\mathbb K(\alpha)$. 
We deduce that $
\prod_{\sigma} \sigma(D) \in \mathbb Z \setminus \{0\}
$ where the product is over all embeddings $\sigma$ of $\mathbb K(\alpha)$ into $\mathbb C$. 

We define for simplicity $\gamma:=f(\alpha)$. Let $L_j:=\sum_{k=0}^\delta a_k \gamma^{k+j}$. On the one hand, by linear combinations of columns, we find  
$$
D=\left\vert 
\begin{matrix}
L_0&a_1&\cdots &a_{\delta}&0&\cdots&\cdots&\cdots &\cdots &0
\\
L_1 &a_0&\cdots &a_{\delta-1}&a_{\delta}&0&\cdots&\cdots&\cdots &0
\\
\vdots &\vdots&\vdots &\vdots&\vdots&\vdots&\vdots&\vdots&\vdots&\vdots
\\
L_{N-\delta} &0&\cdots &0&\cdots&\cdots&0&a_0&\cdots & a_{\delta}
\\
R_{0}&A_{0,1}&\cdots &\cdots &\cdots&\cdots&\cdots &\cdots&\cdots& A_{0,N}
\\
R_{1}&A_{1,1}&\cdots &\cdots &\cdots&\cdots&\cdots &\cdots&\cdots& A_{1,N}
\\
\vdots &\vdots&\vdots &\vdots&\vdots&\vdots&\vdots&\vdots&\vdots&\vdots
\\
R_{\delta-1}&A_{\delta-1,1}&\cdots &\cdots &\cdots& \cdots&\cdots &\cdots &\cdots& A_{\delta-1,N}
\end{matrix}
\right\vert.
$$
Expanding this determinant along its first column and since $L_j=\gamma^j L_0$, the bounds given in \S\ref{sec:ssmethod}  above  
imply that, for all $\varepsilon\in (0,1)$, all $n\ge 0$ and $H \ge 1$, 
\begin{equation} \label{eq:majorationDbis}
\vert D\vert \le C_2 H^{N-\delta}n!^{\delta(1+\varepsilon)} \vert L_0\vert + C_2 \frac{H^{N-\delta+1}}{n!^{N-\delta-1-\delta\varepsilon}}
\end{equation}
for a constant $C_2>0$ independent of $n$ and $H$.

On the other hand, since Eq. \eqref{eqsigmaDun} holds in this setting too, we have  for all $\varepsilon\in (0,1)$, all $n\ge 0$ and $H \ge 1$, 
\begin{equation} \label{eq:majorationsigmaDbis}
\vert \sigma(D)\vert \le C_2H^{N-\delta+1}n!^{\delta(1+\varepsilon)}
\end{equation}
where the constant $C_2>0$ can be taken the same as in Eq. \eqref{eq:majorationDbis} (up to  increasing it if necessary). 

Eqns.  \eqref{eq:majorationDbis} and \eqref{eq:majorationsigmaDbis} are exactly the same as Eqns.  \eqref{eq:majorationD} and \eqref{eq:majorationsigmaD} in \S \ref{sec:determinant}, with $C_2$ instead of $C_1$. The end of the proof of  Theorem~\ref{theo:2} is then exactly the same as the one of Theorem \ref{theo:1}.

\noindent St\'ephane Fischler, Universit\'e Paris-Saclay, CNRS, Laboratoire de math\'ematiques d'Orsay, 91405 Orsay, France.

\medskip

\noindent Tanguy Rivoal, Institut Fourier, Universit\'e Grenoble Alpes, CNRS,  CS 40700, 38058 Grenoble cedex 9, France.

\bigskip

\noindent Keywords: $E$-functions, Transcendence measure, Siegel-Shidlovskii method.

\bigskip

\noindent MSC 2020: 11J82 (Primary), 11J91 (Secondary)

\end{document}